\documentclass{article}
\usepackage{amssymb}  
\usepackage{amsbsy}   
\def\domedskip{\par \ifdim\lastskip<\medskipamount
  \removelastskip \medskip\fi}
\def\qed{\hfill {\hbox{${\vcenter{\vbox{                        
   \hrule height 0.4pt\hbox{\vrule width 0.4pt height 6pt
   \kern5pt\vrule width 0.4pt}\hrule height 0.4pt}}}$}}}
\newenvironment{proof}[1][Proof]{\domedskip\noindent{\bf #1.}\quad}%
{\qed\par\medskip}
\def\sqr#1#2{{\vcenter{\vbox{\hrule height.#2pt
  \hbox{\vrule width.#2pt height#1pt \kern#1pt \vrule width.#2pt}
  \hrule height.#2pt}}}}
\newenvironment{remark}[1][Remark]{\domedskip\noindent{\bf #1.}\quad}%
{\par\domedskip}
\newenvironment{eqnarray+}{\setlength\arraycolsep{.14em}\begin{eqnarray*}}%
{\end{eqnarray*}}
\def\z{\mathbb Z}               
\def\maps#1#2#3{#1\colon #2\to #3} 
\def\id{\mathop{{\rm id}}\nolimits}         
\def\ker{\mathop{{\rm Ker}}\nolimits}       
\def\im{\mathop{{\rm Im}}\nolimits}         
\def\rank{\mathop{{\rm rank}}\nolimits}     
\def\mapright#1{\buildrel #1 \over \longrightarrow}
\def\orb#1{\mathcal{O}_{#1}}
\def\sqset{$\,\sqr65\,$--$\,$set}
\def\tuple#1{\mathbf{#1}}
\def\d{\partial}
\newtheorem{theorem}{Theorem}[section]
\newtheorem{prop}{Proposition}[section]
\newtheorem{lemma}{Lemma}[section]

\title{The Betti numbers of some finite racks}

\author{R. A. Litherland \\
       \textit{\small Department of Mathematics} \\
       \textit{\small Louisiana State University} \\
       \textit{\small Baton Rouge, LA 70803} \\
       \texttt{\small lither@math.lsu.edu} \\
\and 
       Sam Nelson \\
       \textit{\small Department of Mathematics} \\
       \textit{\small Louisiana State University} \\
       \textit{\small Baton Rouge, LA 70803} \\
       \texttt{\small nelson@math.lsu.edu} \\
}

\date{June 18, 2001}

\begin{document}

\maketitle

\begin{abstract}
We show that the lower bounds for Betti numbers given in
\cite{CJKSb} are equalities for a class of racks that includes 
dihedral and Alexander racks. We confirm a conjecture from the
same paper by defining a splitting for the short exact sequence 
of quandle chain complexes. We define isomorphisms between 
Alexander racks of certain forms, and we also list the 
second and third homology groups of some dihedral and Alexander 
quandles.
\end{abstract}

\section{Introduction}

We start by recalling some basic definitions. 
Let $X$ be a non-empty set with a binary operation,
which, following Fenn and Rourke \cite{FR}, we write as 
exponentiation: $(a,b) \mapsto a^b$. This allows us to dispense with brackets
by using the standard conventions that $a^{b^c}$ means $a^{(b^c)}$
and $a^{bc}$ means $(a^b)^c$. Then $X$ is a \textit{rack} if it satisfies the
following two axioms.
\begin{itemize}
  \item[(i)] For all $a, b \in X$, there is a unique element $c$ of $X$
such that $c^a = b$. 
 \item[(ii)] (The rack identity) For all $a, b, c \in X$, $a^{bc} = a^{cb^c}$.
\end{itemize}
A \textit{quandle} is a rack satisfying one further axiom.
\begin{itemize}
\item[(iii)] (The quandle condition) For all $a \in X$, 
$a^a = a$.
\end{itemize}
A rack is \textit{trivial} if $a^b = a$ for all $a$ and $b$.

By axiom (i), the function $\maps{f_a}XX$ defined by $f_a(b) = b^a$
is a bijection.
For $a, b \in X$, we set $a^{\bar b}= f_b^{-1}(a)$. 
Here $\bar b$ does not denote
an element of $X$, but we may identify $\bar b$ with the inverse of $b$
in the free group $F(X)$ on $X$. This allows us to define a (right) action
of $F(X)$ on $X$, and by an \textit{orbit} of $X$ we mean an orbit under this 
action. The set of orbits of $X$ will be denoted by $\orb X$, and the 
projection from $X$ to $\orb X$ by $\pi$. We regard $\orb X$ as a trivial
rack, and then $\pi$ is a rack homomorphism.

We now define the class of racks that we shall study in \S3 of this paper.
Let $X$ be a finite rack, and $a, b \in X$. Let $N(a,b)$ be the
number of elements $c$ of $X$ such that $a^c = b$. Of course,
$N(a, b) = 0$ if $a$ and $b$ are in different orbits. We say that
$X$ \textit{has homogeneous orbits} if, for each orbit $\omega$ and
each pair of elements $a$ and $b$ of $\omega$, $N(a,b)$ depends only
on $\omega$. If this is so, then $|\omega|$ divides $|X|$ for each
$\omega \in \orb X$, and $N(a,b) = |X|/|\omega|$ for all
$a, b \in \omega$; we set $N_\omega = |X|/|\omega|$.

Let us consider some of the standard examples of racks in the light of
this definition. Clearly (if uninterestingly), 
any trivial finite rack has
homogeneous orbits. So does any finite conjugation rack $\hbox{conj}(G)$.
(Here $G$ is a group, and $\hbox{conj}(G)$ denotes $G$ with
the rack operation $g^h = h^{-1}gh$.) Fenn and Rourke use the term
\textit{conjugation rack} in a broader sense, to refer to any union of
conjugacy classes in a group. In general, these do not have homogeneous
orbits (consider $G-\{1\}$); 
however, any dihedral rack $R_n$ does. ($R_n$ is the set of
reflections in the dihedral group of order $2n$.) This is easy to verify
directly, and also follows from Proposition 1 below. Any  cyclic
rack (except the trivial rack of order 1) does not have homogeneous orbits. 
(The cyclic rack $C_n$ of order $n$
is the set $\{\,0, 1, \ldots, n-1\,\}$ with the operation
$a^b = a+1\bmod n$. Here there is only one orbit, but $N(a,b) = n$
if $b = a+1\bmod n$, and is 0 otherwise.)

As an example of a non-quandle that does have homogeneous orbits,
consider a four-element set $X=\{\,a,b,c,d\,\}$. We define the operation
by specifying the permutations $f_x$ of $X$: $f_a=f_b$ is the transposition
exchanging $a$ and $b$, and $f_c=f_d$ is the identity. One may check
that the rack identity holds, most easily by using the third form of the 
identity given in \cite{FR}; the quandle condition clearly does not.
The only non-trivial orbit is $\{\,a,b\,\}$, and
$N(a,a) = N(a,b) = N(b,a) = N(b,b) = 2$.

Next we consider the finite Alexander racks. Let $M$ be any module
over the ring $\Lambda = \z[t, t^{-1}]$ of one-variable Laurent polynomials.
Then $M$ may be made into a rack by the operation $a^b = ta+(1-t)b$, and a
rack obtained this way is called an \textit{Alexander rack}. For
$M=\z_n[t,t^{-1}]/(t+1)$, the Alexander rack is isomorphic to $R_n$.

\begin{prop}\label{prop:alex} Let $M$ be a finite $\Lambda$-module, and 
let $\bar M$ be the quotient of $M$ by the submodule $(1-t)M$. 
When $M$ is considered as an Alexander rack:
\begin{itemize}
\item[\rm (a)] $M$ has homogeneous orbits; and
\item[\rm (b)] $\orb M$ may be identified with $\bar M$.
\end{itemize}
\end{prop}
\begin{proof} 
Let $\maps p M {\bar M}$ be the natural map. We have $a^x = a^y$ iff
$(1-t)(x-y)=0$, so for any $a, b\in M$, $N(a,b)$ is either 0 or the order
of $\ker(\maps {1-t} M M)$. The result will follow once we show that,
for $a,b\in M$, the following statements are equivalent:
\begin{itemize}
\item[(1)] $a$ and $b$ are in the same orbit;
\item[(2)] $p(a)=p(b)$;
\item[(3)] $N(a,b) \neq 0$.
\end{itemize}
Now $a-a^b = (1-t)(a-b)$, so $p(a) = p(a^b)$, from which it follows that
(1) implies (2). If $p(a) = p(b)$, then $b = a + (1-t)c$ for some
$c \in M$, which gives $b = a^{a+c}$. Thus (2) implies (3), and trivially
(3) implies (1).
\end{proof}

In ~\cite{FRS}, Fenn, Rourke and Sanderson associate to each rack $X$ a
\sqset\ (a cubical set without degeneracies) as follows.
The set of $n$-cubes is $X^n$, and the face maps 
$\maps {\d^\epsilon_i}{X^n}{X^{n-1}}$ ($1 \leq i \leq n$, 
$\epsilon = 0$ or 1) are defined by
\begin{eqnarray+}
 \d^0_i(x_1,\ldots,x_n) & = & (x_1,\ldots,x_{i-1},x_{i+1},\ldots,x_n); \\
 \d^1_i(x_1,\ldots,x_n) & = & (x_1^{x_i},\ldots,x_{i-1}^{x_i},x_{i+1},\ldots,x_n). 
\end{eqnarray+}

We follow Carter, Jelsovsky, Kamada and Saito \cite{CJKSb} in denoting
the associated chain complex by $C^R_*(X)$, and calling its homology
$H^R_*(X)$ the \textit{rack homology} of $X$.  Thus $C^R_n(X)$ is the
free abelian group on $X^n$, and the boundary map
$\maps\d{C^R_n(X)}{C^R_{n-1}(X)}$ is defined by $\d =
\sum_{i=1}^n(-1)^i(\d_i^0 - \d_i^1)$.  Now suppose that $X$ is a
quandle, and define $C^D_n(X)$ to be the subgroup of $C^R_n(X)$
generated by $n$-tuples $(x_1,\ldots,x_n)$ with $x_i = x_{i+1}$ for
some $i$, $1 \leq i < n$. It follows from the quandle condition that
$C^D_*(X)$ is a subcomplex of $C^R_*(X)$. The quotient complex is
denoted by $C^Q_*(X)$, and its homology $H^Q_*(X)$ is called the
\textit{quandle homology} of $X$. The homology $H^D_*(X)$ of
$C^D_*(X)$ is the {\it degeneration homology\/} of $X$.  We shall use
the convention that in an expression such as $C^W_n(X)$, $W$ may be
any one of $R$, $Q$ or $D$ if $X$ is a quandle, but is always $R$ if
not.  There are Betti numbers $\beta^W_n(X) = \rank H^W_n(X)$. There
are also homology and cohomology groups with coefficients in any
abelian group $G$, denoted by $H^W_n(X;G)$ and $H_W^n(X;G)$. For the
applications to knot theory, the groups of interest are the cohomology
groups with coefficients in $\z_p$ (the integers modulo a prime $p$),
but since these are determined by the integral homology groups we
shall concentrate on the latter.  The homology groups in dimensions 0
and 1 are easily computed; see Proposition 3.8 of \cite{CJKSb}.  
When the set of orbits of $X$ is regarded as a trivial rack, the chain complex
$C^W_*(\orb X)$ has all its boundary maps zero, so $H^W_n(\orb X) =
C^W_n(X)$. Thus when $X$ is a finite rack with $m$ orbits, $H^W_n(\orb X)$
($n \geq 1$)
is a free abelian group of rank $m^n$, $m(m-1)^{n-1}$ or
$m^n-m(m-1)^{n-1}$ for $W=R$, $Q$ or $D$, respectively. In
\cite{CJKSb}, it is shown that in this case $\beta^W_n(X) \geq
\beta^W_n(\orb X)$. 
(It is not explicitly stated in \cite{CJKSb} that the case $W=R$ holds
when $X$ is not a quandle, but this is so
by essentially the same proof.) We now state our main result, which
shows that these bounds are exact in many cases.

\begin{theorem} \label{thm:main}
Let $X$ be a finite rack with homogeneous orbits.
Then $\beta^W_n(X) = \beta^W_n(\orb X)$,
and the torsion subgroup of $H^W_n(X)$ is annihilated by $|X|^n$.
\end{theorem}

\begin{remark} While this paper was in preparation, we learned that
Mochizuki has proved an almost identical theorem by a different method
(\cite{M}, Theorem 1.1). The main difference in the results is that
Mochizuki's theorem applies only to finite Alexander racks.
\end{remark}

The case $W=R$ of Theorem 1.1 is proved directly. For the other cases, we
need to prove conjecture 3.11 of \cite{CJKSb}; this is done in
\S\ref{sec:split}. Theorem \ref{thm:main} is proved in
\S\ref{sec:main}, and in \S\ref{sec:calc} we report on some machine
calculations of homology groups.

\section{Splitting the difference between quandle and rack homology}
\label{sec:split}

In this section, $X$ will always denote a quandle. Also, we
redefine $C^R_0(X)$ and $C^Q_0(X)$ to be 0. ($C^D_0(X)$ is already 0.) 
This loses no information, and allows us to avoid treating dimension 0 
as a special case at various points. Strictly speaking, we shall be 
working with the reduced complexes $\tilde C^R_*(X)$ and 
$\tilde C^Q_*(X)$, but we abuse notation by leaving off the tildes. 
From the short exact sequence
\begin{equation}
 0 \to C^D_*(X) \to C^R_*(X) \to C^Q_*(X) \to 0 
\end{equation} 
of chain complexes, we have a long exact sequence
\[ 
 \dots \to H^D_n(X) \to H^R_n(X) \to H^Q_n(X) \to H^D_{n-1}(X) \to \dots
\] 
of homology groups. In \cite{CJKSb} it is proved (in
Proposition 3.9) that the connecting homomorphism
$H^Q_n(X)\to H^D_{n-1}(X)$ is the zero map when $n = 3$, and
conjectured that this is so for all $n$; in \cite{CJKSd} (Theorem 8.2) the
case $n=4$ is proved. We show that the conjecture is indeed true; in
fact we prove more.

\begin{theorem} \label{thm:split}
For any quandle $X$, the short exact sequence 
{\rm (1)} is split.
\end{theorem}
\begin{remark} It is easy to see that, for each $n$, the sequence
\[
 0 \to C^D_n(X) \to C^R_n(X) \to C^Q_n(X) \to 0 
\]
of abelian groups is split, but the obvious splittings are not
compatible with the boundary maps.
\end{remark}

If $\tuple x = (x_1, \ldots, x_n)\in X^n$ and $y \in X$, we set
$\tuple x * y = (x_1, \ldots, x_n, y)\in X^{n+1}$ and $\tuple x^y =
(x_1^y, \ldots, x_n^y) \in X^n$. Then, for $c \in C^R_n(X)$ we define
$c * y \in C^R_{n+1}(X)$ and $c^y\in C^R_n(X)$ by linearity in
$c$. Note that $\d(c * y) = \d(c)*y +(-1)^{n+1}(c - c^y)$.
Next we define homomorphisms $\maps{\alpha_n}{C^R_n(X)}{C^R_n(X)}$ by
induction on $n$. We take $\alpha_1$ to be the identity map, and for
$n \geq 1$, $\tuple x \in X^n$ and $y\in X$ we set
\[
 \alpha_{n+1}(\tuple x * y) = \alpha_n(\tuple x)*y -\alpha_n(\tuple x)*x_n.
\]
We also define homomorphisms $\maps{\beta_n}{C^R_n(X)}{C^R_{n+1}(X)}$
by $\beta_n(\tuple x) = \alpha_n(\tuple x)*x_n$.  Then, for any $c \in
C^R_n(X)$ and $y\in X$ we have
\[
 \alpha_{n+1}(c * y) = \alpha_n(c)*y - \beta_n(c).
\]

\begin{lemma} The homomorphisms $\maps{\alpha_n}{C^R_n(X)}{C^R_n(X)}$
form a chain map $\maps{\alpha}{C^R_*(X)}{C^R_*(X)}$.
\end{lemma}
\begin{proof} Note first that for $\tuple x \in X^n$ and $y\in X$ we have
$\alpha_n(\tuple x^y) = \alpha_n(\tuple x)^y$. We prove that
$\d\alpha_n = \alpha_{n-1}\d$ by induction on $n \geq 2$.
For $n = 2$ we have $\alpha_2(x,y) = (x,y)-(x,x)$, so since $(x,x)$ is
a cycle, $\d\alpha_2(x,y) = \d(x,y) =
\alpha_1\d(x,y)$.  Suppose then that the result is true for some
$n \geq 2$, and let $\tuple x \in X^n$ and $y\in X$. We compute
\begin{eqnarray+}
\d\alpha_{n+1}(\tuple x * y)
&=&\d\bigl(\alpha_n(\tuple x) * y\bigr) -
\d\bigl(\alpha_n(\tuple x) * x_n\bigr)\\ 
&=&\d\alpha_n(\tuple x) * y + (-1)^{n+1}\bigl(\alpha_n(\tuple x) -
\alpha_n(\tuple x)^y\bigr) \\
&&\quad{}- \d\alpha_n(\tuple x) * x_n 
- (-1)^{n+1}\bigl(\alpha_n(\tuple x) 
- \alpha_n(\tuple x)^{x_n}\bigr) \\
&=&\alpha_{n-1}\d(\tuple x) * y -
\alpha_{n-1}\d(\tuple x) * x_n 
+ (-1)^n\alpha_n(\tuple x^y) -
(-1)^n\alpha_n(\tuple x^{x_n})
\end{eqnarray+}
and
\begin{eqnarray+}
\alpha_n\d(\tuple x * y) &=& \alpha_n\bigl(\d(\tuple x)*y 
+ (-1)^{n+1}(\tuple x - \tuple x^y)\bigr)\\ 
&=& \alpha_{n-1}\d(\tuple x)*y - \beta_{n-1}\d(\tuple x) 
- (-1)^n\alpha_n(\tuple x) + (-1)^n\alpha_n(\tuple x^y).
\end{eqnarray+} 
Hence $\d\alpha_{n+1}(\tuple x *y) = \alpha_n\d(\tuple x * y)$ iff
$$\alpha_{n-1}\d(\tuple x) * x_n + (-1)^n\alpha_n(\tuple x^{x_n}) 
= \beta_{n-1}\d(\tuple x) + (-1)^n\alpha_n(\tuple x).\eqno{(2)}$$ 
Now, for $1 \leq i < n$ and $\epsilon = 0$ or 1,
$\d_i^\epsilon(\tuple x)$ is an element of $X^{n-1}$ with last
entry $x_n$, so $\alpha_{n-1}\d_i^\epsilon(\tuple x)*x_n
=\beta_{n-1}\d_i^\epsilon(\tuple x)$. Further,
\begin{eqnarray+}
\alpha_{n-1}\d_n^0(\tuple x)*x_n-\beta_{n-1}\d_n^0(\tuple x) 
&=& \alpha_n\bigl(\d_n^0(\tuple x)*x_n\bigr) =\alpha_n(\tuple x)\\
\hbox{and}\qquad 
\alpha_{n-1}\d_n^1(\tuple x)*x_n -\beta_{n-1}\d_n^1(\tuple x) 
&=&\alpha_n\bigl(\d_n^1(\tuple x)*x_n\bigr)
=\alpha_n(\tuple x^{x_n}).
\end{eqnarray+}
(The last step here uses the quandle condition.)  It
follows that
\[ \alpha_{n-1}\d(\tuple x)*x_n -\beta_{n-1}\d(\tuple x)
 =(-1)^n\bigl(\alpha_n(\tuple x)-\alpha_n(\tuple x^{x_n})\bigr),
\]
proving equation (2), and with it the lemma.
\end{proof}

\begin{proof}[Proof of Theorem \ref{thm:split}] We show that the chain map
$C^R_*(X) \to C^R_*(X)$ sending $c$ to $c - \alpha(c)$ is a projection
onto the subcomplex $C^D_*(X)$. We must prove the following two
statements.
\begin{itemize}
\item[(a)] If $c \in C^D_n(X)$ then $\alpha_n(c) = 0$.
\item[(b)] If $c \in C^R_n(X)$ then $c-\alpha_n(c) \in C^D_n(X)$.
\end{itemize}
For $n=1$, $C^D_n(X)=0$, so (a) is true in this case.  
Let $\tuple x\in X^n$ ($n \geq 1$) and $y\in X$, and suppose that 
$\tuple x *y \in C^D_{n+1}(X)$.  Then either $\tuple x \in C^D_n(X)$
or $x_n = y$, and
it follows that $\alpha_{n+1}(\tuple x * y) = 0$ (using induction in
the first case).  Thus (a) is proved.  As for (b), this is clear for
$n=1$, so suppose that it holds for some $n \geq 1$ and take 
$\tuple x \in X^n$ and $y\in X$. Then
\begin{eqnarray+}
\tuple x * y - \alpha_{n+1}(\tuple x * y) - \tuple x * x_n
&=& \bigl(\tuple x -\alpha_n(\tuple x)\bigr)*y 
-\bigl(\tuple x -\alpha_n(\tuple x)\bigr)*x_n\\ 
&\in &C^D_{n+1}(X).
\end{eqnarray+} 
Since $\tuple x * x_n$ is in $C^D_{n+1}(X)$, so is $\tuple x * y -
\alpha_{n+1}(\tuple x * y)$, and (b) follows.
\end{proof}
We shall denote the free abelian group on a set $A$ by $\z[A]$.
(This is consistent with the usage $\z[G]$ for a group ring.)
It is shown in Proposition 3.9 of \cite{CJKSb} that $H^D_2(X) \cong
\z[\orb X]$.  Combining this with Theorem \ref{thm:split} 
gives the first assertion
of the next theorem; for the second we need some lemmas.

\begin{theorem}\label{thm:h23} For any quandle $X$, we have
\begin{eqnarray+}
H^R_2(X) &\cong & H^Q_2(X) \oplus \z[\orb X]\\
\hbox{and}\qquad 
H^R_3(X) &\cong & H^Q_3(X) \oplus H^Q_2(X) \oplus \z[\orb X^2].\\
\end{eqnarray+}
\end{theorem}

Let $C^L_n(X)$ be the subgroup of $C^D_n(X)$ generated by $n$-tuples
$(x_1, \ldots ,x_n)$ with $x_i=x_{i+1}$ for some $i$, $2\leq i <
n$. (We use the letter $L$ because the degeneracy occurs late in these
$n$-tuples.) Note that $C^L_n(X) = 0$ for $n < 3$.

\begin{lemma} The subgroups $C^L_n(X)$ form a subcomplex of
$C^D_*(X)$.
\end{lemma}
\begin{proof} Let $\tuple x = (x_1, \ldots ,x_n)$ have $x_i=x_{i+1}$ for 
some $i$ with $2\leq i < n$.  Since $\d_i^\epsilon(\tuple x) =
\d_{i+1}^\epsilon(\tuple x)$ for $\epsilon = 0$ or 1 (and, as
for any $\tuple x \in X^n$, $\d_1^0(\tuple x) =
\d_1^1(\tuple x)$), we have
$$\d(\tuple x) =
\sum_{j=2}^{i-1}(-1)^j\bigl(\d_j^0(\tuple x) 
- \d_j^1(\tuple x)\bigr) +
\sum_{j=i+2}^n(-1)^j\bigl(\d_j^0(\tuple x) 
- \d_j^1(\tuple x)\bigr).\eqno{(3)}$$ 
Fix $j$ and $\epsilon$, and set $\tuple y =
(y_1,\ldots,y_{n-1}) = \d_j^\epsilon(\tuple x)$.  If $i = 2$,
the first sum in (3) is empty. If $i > 2$ and $2 \leq j \leq i-1$,
$y_{i-1} = y_i$, so $\tuple y\in C^L_{n-1}(X)$.  For $i+2 \leq j \leq
n$, $y_i = y_{i+1}$, so again $\tuple y\in C^L_{n-1}(X)$, and it
follows that $\d(\tuple x)\in C^L_{n-1}(X)$.
\end{proof}

\begin{lemma}\label{lem:cqdl} There is an isomorphism of chain complexes
$C^D_*(X) \cong C^Q_{*-1}(X) \oplus C^L_*(X)$.
\end{lemma}

\begin{proof} 
We let $\maps{i}{C^D_*(X)}{C^R_*(X)}$ and
$\maps{j}{C^L_*(X)}{C^D_*(X)}$ be the inclusions.  Define
$\maps{r}{C^R_{*-1}(X)}{C^D_*(X)}$ by 
\[r_n(x_1,x_2,\ldots,x_{n-1}) =
(x_1,x_1,x_2,\ldots,x_{n-1})\]
for $n \geq 2$.  (For $n \leq 1$ the
groups involved are 0.) A straightforward computation shows that $r$
is a chain map.  Since $r(C^D_{*-1}(X)) \leq C^L_*(X)$, $r$ induces
$\maps{s}{C^D_{*-1}(X)}{C^L_*(X)}$.

Now $r$ is injective, $C^D_*(X)$ is generated by $\im(r)$ and
$C^L_*(X)$, and $\im(r)\cap C^L_*(X) = \im(r\circ i) = \im(j \circ
s)$. Hence there is a short exact sequence
$$0 \to C^D_{*-1}(X) \mapright{\phi} C^R_{*-1}(X) \oplus C^L_*(X)
\mapright{\psi} C^D_*(X) \to 0,$$
where $\phi(c) = (i(c), -s(c))$ and $\psi(d,e) = r(d) + j(e)$. By
Theorem 1, there is an isomorphism of chain complexes
$\maps{\chi}{C^R_{*-1}(X)}{C^Q_{*-1}(X)
\oplus C^D_{*-1}(X)}$ such that, for $c \in C^D_{*-1}(X)$,
$\chi i(c) = (0, c)$. Then $C^D_*(X)$ is isomorphic to the cokernel of
$$\maps{(\chi\oplus \id)\circ\phi}{C^D_{*-1}(X)}{C^Q_{*-1}(X)
\oplus C^D_{*-1}(X)\oplus C^L_*(X)}.$$
But, for $c \in C^D_{*-1}(X)$, 
$(\chi\oplus \id)(\phi(c)) = (0, c, -s(c))$, 
so this cokernel is isomorphic to $C^Q_{*-1}(X) \oplus C^L_*(X)$, 
and we are done.
\end{proof}

We denote the homology of $C^L_*(X)$ by $H^L_*(X)$.
\begin{lemma} \label{lem:h3l}
For any quandle $X$, $H^L_3(X) \cong \z[\orb X^2]$.
\end{lemma}
\begin{proof} A basis for $C^L_3(X)$ consists of all elements of $X^3$ of the form
$(x, y, y)$, and these are all cycles. The group $C^L_4(X)$ is
generated by all elements of $X^4$ of one of the forms $(x,y,y,z)$ and
$(x,z,y,y)$, and we have
\begin{eqnarray+}
\d(x,y,y,z) &=& (x,y,y) - (x^z,y^z,y^z)\\
\hbox{and}\qquad \d(x,z,y,y) &=& (x,y,y) - (x^z,y,y).
\end{eqnarray+}
It follows that $H^L_3(X)$ is free abelian, with a basis consisting of
the equivalence classes of triples $(x,y,y)$ under the equivalence
relation $\sim$ generated by
\[ (x,y,y) \sim (x^z, y, y) \sim(x^z, y^z, y^z)\qquad
 \hbox{for all $x,y,z\in X$.}
\]
Given $x,y,z \in X$, let $w$ be the element of $X$ such that $w^z =
x$.  Then $(w,y,y)\sim(w^z,y,y) = (x,y,y)$ and $(w,y,y)\sim (w^z,
y^z,y^z) =(x,y^z,y^z)$, so $(x,y,y)\sim(x,y^z,y^z)$. It follows that
$(x,y,y)\sim (x',y',y')$ iff $\pi(x)=\pi(x')$ and $\pi(y)=\pi(y')$, so
the set of equivalence classes of $\sim$ may be identified with 
$\orb X^2$.
\end{proof}

The second assertion of Theorem \ref{thm:h23} follows immediately from 
Theorem \ref{thm:split} and Lemmas \ref{lem:cqdl} and \ref{lem:h3l}.

\section{Proof of Theorem \ref{thm:main}}
\label{sec:main}

In this section, $X$ is a rack with homogeneous orbits, and 
$\tuple x = (x_1,\ldots,x_n)$ is an element of $X^n$ ($n >= 0$).
Define $\phi_n^j:C^R_n(X)\to C^R_n(X)$ by
\[\phi_n^j(\tuple{x}) =\left\{ \begin{array}{ll}
                  \tuple{x}     
            & \mbox{for $j=0$} \\
                  \sum_{\tuple{y}\in X^j} (x_1^{y_1},\dots,x_j^{y_j},x_{j+1},\dots, x_n) & \mbox{for $1\le j\le n$}\\
                  |X|^{j-n}\phi_n^n(\tuple{x})
            & \mbox{for $j > n$}
                   \end{array}
           \right. \]
and  $D_n^j:C^R_n(X)\to C^R_{n+1}(X)$ by
\[D_n^j(\tuple{x}) =\left\{ \begin{array}{ll}
            \sum_{\tuple{y}\in X^j}
                  (x_1^{y_1},\dots,x_{j-1}^{y_{j-1}},x_j, y_j,
                  x_{j+1},\dots, x_n) & \mbox{for $1\le j\le n$}\\
                  0
            & \mbox{for $j > n$.}
                   \end{array}
           \right. \]
Note that $D_n^1 = \sum_{y\in X} (x_1,y,x_2,\dots,x_n).$

We have homomorphisms of graded groups $\maps{\phi^j =
(\phi^j_n)_{n=0}^\infty}{C^R_*(X)}{C^R_*(X)}$ for $j \ge 0$
and $\maps{D^j =(D^j_n)_{n=0}^\infty}{C^R_*(X)}{C^R_{*+1}(X)}$
for $j \ge 1$.
We will show though a series of lemmas that $D^j$ is a chain homotopy
carrying $\phi^j$ to $|X|\phi^{j-1}$, and hence each $\phi^j$ is
chain homotopic to $|X|^j$ times the identity. Note that this also implies 
$\phi^j$ is a chain map. 

\begin{lemma} \label{lem:sumL}
Let $G$ be an abelian group. Then if $g\colon X\to G$ is
a function we have
\[\sum_{y\in X}g(x^y) = \sum_{y\in X} g(x^{yw})\] 
for any word $w\in F(X)$ in the free group on $X$.
\end{lemma}

\begin{proof} As $y$ runs over $X$, $x^y$ runs over $\pi(x)$, taking on each value 
$N_{\pi(x)}$ times. Thus 
\[ \sum_{y\in X} g(x^y)=N_{\pi(x)}\sum_{z\in \pi(x)} g(z).\]
The automorphism $f_w:X\to X$ given by $f_w(x)=x^w$ is in particular 
a bijection and carries $\pi(x)$ to itself, so the restriction 
$f\mid_{\pi(x)}$ is also a bijection. Hence the sum
\[ \sum_{y\in X}g(x^{yw}) =  \sum_{y\in X} g(f_w(x^y))=
N_{\pi(x)}\sum_{z\in \pi(x)} g(z) = \sum_{y\in X} g(x^y). \]
\end{proof}

\begin{lemma} \label{lem:sumR}
Let $G$ be an abelian group. Then if $g\colon X\to G$ 
is a function we have
\[\sum_{y\in X}g(x^y) = \sum_{y\in X} g(x^{wy})\] 
for any word $w\in F(X)$ in the free group on $X$.
\end{lemma}

\begin{proof} Since $\pi(x^w)=\pi(x)$, we have
\[ \sum_{y\in X}g(x^y) 
= N_{\pi(x)}\left( \sum_{z\in \pi(x)}g(z)\right) = N_{\pi(x^w)}  
 \left( \sum_{z\in \pi(x^w)}g(z)\right) = \sum_{y\in X}g(x^{wy}). \]
\end{proof}

\begin{lemma} \label{lem:dDlow}
For $1\le i\le j \le n$, $\d_i^0D_n^j(\tuple{x})=
 \d_i^1D_n^j(\tuple{x})$.
\end{lemma}

\begin{proof} For $i<j$, we have
\begin{eqnarray+}
 \d_i^0D_n^j(\tuple{x}) & = & \d_i^0\left( \sum_{\tuple{y}\in X^j} 
 (x_1^{y_1},\dots, x_{j-1}^{y_{j-1}},x_j, y_j, x_{j+1},\dots, x_n) \right) \\
 & = & \sum_{\tuple{y}\in X^j}(x_1^{y_1},\dots,x_{i-1}^{y_{i-1}},
 x_{i+1}^{y_{i+1}},\dots,x_{j-1}^{y_{j-1}},x_j,y_j,\dots,x_n)\end{eqnarray+}
and
\begin{eqnarray+}
 \d_i^1D_n^j(\tuple{x}) & = & \d_i^1 \left( \sum_{\tuple{y}\in X^j}  
 (x_1^{y_1},\dots, x_{j-1}^{y_{j-1}},x_j, y_j, x_{j+1},\dots, x_n) \right) \\ 
 & = & \sum_{\tuple{y}\in X^j}(x_1^{y_1x_i^{y_i}},\dots,
 x_{i-1}^{y_{i-1}x_i^{y_i}},x_{i+1}^{y_{i+1}},\dots,x_{j-1}^{y_{j-1}},
 x_j,y_j,\dots,x_n).\end{eqnarray+}
For $i=j$ we have
\[ \d_j^0D_n^j(\tuple{x}) = \sum_{\tuple{y}\in X^j} (x_1^{y_1},\dots,
 x_{j-1}^{y_{j-1}},y_j,\dots,x_n) \]
and 
\[ \d_j^1D_n^j(\tuple{x}) = \sum_{\tuple{y}\in X^j} (x_1^{y_1x_j},\dots,
 x_{j-1}^{y_{j-1}x_j},y_j,\dots,x_n).\]
Applying Lemma \ref{lem:sumL} $i-1$ times, the sums agree as required.
\end{proof}

\begin{lemma} \label{lem:Ddlow}
For $1 \le i\le j < n$, $D_{n-1}^j\d_i^0(\tuple{x})=
D_{n-1}^j\d_i^1(\tuple{x})$.
\end{lemma}

\begin{proof} For $1\le i\le j$, 
\begin{eqnarray+} 
 D_{n-1}^j\d_i^0(\tuple{x}) & = & D_{n-1}^j (x_1,\dots,x_{i-1},
 x_{i+1},\dots,x_n)\\
 & = & \sum_{\tuple{y}\in X^j}(x_1^{y_1},\dots,x_{i-1}^{y_{i-1}},
 x_{i+1}^{y_i}, \dots,x_j^{y_{j-1}}, x_{j+1},y_j,x_{j+2},\dots,x_n) 
\end{eqnarray+}
and
\begin{eqnarray+} 
 D_{n-1}^j\d_i^1(\tuple{x}) & = & D_{n-1}^j (x_1^{y_1},\dots,
 x_{i-1}^{y_{i-1}}, x_{i+1},\dots,x_n)\\
 & = & \sum_{\tuple{y}\in X^j}(x_1^{x_iy_1},\dots,x_{i-1}^{x_iy_{i-1}},
 x_{i+1}^{y_i}, \dots,x_j^{y_{j-1}}, x_{j+1},y_j,x_{j+2},\dots,x_n). 
\end{eqnarray+}

Applying Lemma \ref{lem:sumR} $i-1$ times, the sums agree as required.
\end{proof}

\begin{lemma} \label{lem:dDmid0}
For $1\le j \le n$, 
$\d_{j+1}^0D_n^j(\tuple{x})=|X|\phi_n^{j-1}(\tuple{x})$.
\end{lemma}

\begin{proof}
\begin{eqnarray+} 
 \d_{j+1}^0D_n^j(\tuple{x}) & = & \d_{j+1}^0\left( \sum_{\tuple{y}\in X^j} 
 (x_1^{y_1},\dots,x_{j-1}^{y_{j-1}},x_j, y_j, x_{j+1},\dots, x_n)\right) \\
 & = & \sum_{\tuple{y}\in X^j} (x_1^{y_1},\dots,x_{j-1}^{y_{j-1}},x_j, 
 x_{j+1},\dots, x_n) \\
 & = & \sum_{y_j\in X} \phi_n^{j-1}(\tuple{x}) \\ 
 & = & |X|\phi_n^{j-1}(\tuple{x}).\  
\end{eqnarray+}
\end{proof}

\begin{lemma} \label{lem:dDmid1}
For $1\le j \le n$, 
$\d_{j+1}^1D_n^j(\tuple{x})=\phi_n^{j}(\tuple{x})$.
\end{lemma}

\begin{proof}
\begin{eqnarray+} 
 \d_{j+1}^1D_n^j(\tuple{x}) & = & \d_{j+1}^1\left( \sum_{\tuple{y}\in X^j} 
 (x_1^{y_1},\dots,x_{j-1}^{y_{j-1}},x_j, y_j, x_{j+1},\dots, x_n)\right) \\
 & = & \sum_{\tuple{y}\in X^j} (x_1^{y_1y_j},\dots,x_{j-1}^{y_{j-1}y_j},
 x_j^{y_j}, x_{j+1},\dots, x_n) \\
 & = & \sum_{\tuple{y}\in X^j} (x_1^{y_1},\dots, x_{j-1}^{y_{j-1}},x_j^{y_j}, 
 x_{j+1},\dots, x_n) \\
 & = & \phi_n^{j}(\tuple{x})
\end{eqnarray+}
by $j-1$ applications of Lemma \ref{lem:sumL}.
\end{proof}

\begin{lemma} \label{lem:dDhigh0}
For $1 \le j < i \le n$, $D_{n-1}^j\d_i^0(\tuple{x})=
 \d_{i+1}^0D_n^j(\tuple{x}).$
\end{lemma}

\begin{proof}
\begin{eqnarray+} 
 D_{n-1}^j\d_i^0(\tuple{x}) & = & D_{n-1}^j (x_1,\dots,x_{i-1}, 
 x_{i+1},\dots,x_n) \\
 & = & \sum_{\tuple{y}\in X^j}(x_1^{y_1},\dots,x_{j-1}^{y_{j-1}},x_j,
 y_j,x_{j+1},\dots,x_{i-1}, x_{i+1}, \dots,x_n) 
\end{eqnarray+} and
\begin{eqnarray+} 
 \d_{i+1}^0D_n^j(\tuple{x}) & = & \d_{i+1}^0\left( \sum_{\tuple{y}\in X^j} 
 (x_1^{y_1},\dots,x_{j-1}^{y_{j-1}},x_j, y_j, x_{j+1},\dots, x_n) \right) \\
 & = & \sum_{\tuple{y}\in X^j} (x_1^{y_1},\dots,x_{j-1}^{y_{j-1}},x_j, y_j, 
 x_{j+1},\dots,x_{i-1}, x_{i+1},\dots, x_n). 
\end{eqnarray+}
\end{proof}

\begin{lemma} \label{lem:dDhigh1}
For $1 \le j < i \le n$, $D_{n-1}^j\d_i^1(\tuple{x})=
 \d_{i+1}^1D_n^j(\tuple{x})$.
\end{lemma}

\begin{proof}
\begin{eqnarray+} 
 D_{n-1}^j\d_i^1(\tuple{x}) & = &
 D_{n-1}^j (x_1^{x_i},\dots, x_{i-1}^{x_i}, x_{i+1},\dots,x_n) \\
 & = & \sum_{\tuple{y}\in X^j} (x_1^{x_iy_1},\dots,x_{j-1}^{x_iy_{j-1}},
 x_j^{x_i},y_j,x_{j+1}^{x_i},\dots,x_{i-1}^{x_i},x_{i+1},\dots,x_n) \\ 
 & = & \sum_{\tuple{y}\in X^j} (x_1^{y_1},\dots,x_{j-1}^{y_{j-1}},
 x_j^{x_i},y_j,x_{j+1}^{x_i},\dots,x_{i-1}^{x_i},x_{i+1},\dots,x_n) 
\end{eqnarray+} 
and
\begin{eqnarray+} 
 \d_{i+1}^1D_n^j(\tuple{x}) & = & \d_{i+1}^1 \left(\sum_{\tuple{y}\in X^j} 
 (x_1^{y_1},\dots,x_{j-1}^{y_{j-1}},x_j, y_j, x_{j+1},\dots, x_n) \right) \\
 & = & \sum_{\tuple{y}\in X^j} (x_1^{y_1x_i},\dots,x_{j-1}^{y_{j-1}x_i},
 x_j^{x_i},y_j^{x_i},x_{j+1}^{x_i},\dots,x_{i-1}^{x_i},x_{i+1},\dots,x_n), \\
 & = & \sum_{\tuple{y}\in X^j} (x_1^{y_1},\dots,x_{j-1}^{y_{j-1}},
 x_j^{x_i},y_j^{x_i},x_{j+1}^{x_i},\dots,x_{i-1}^{x_i},x_{i+1},\dots,x_n) 
\end{eqnarray+}
by $j-1$ applications of Lemmas \ref{lem:sumL} and \ref{lem:sumR}. 
But these sums agree as
the set $\{\,y_j^{x_i} \mid y_j\in X\,\}$ is the image 
of $\{\,y_j \mid y_j\in X\,\}$
under the bijection $f_{x_i}$. 
\end{proof}

Putting all this together, we have

\begin{prop}\label{prop:htic} For $j \ge 1$, 
$D^j:C_*^R(X)\to C_{*+1}^R(X)$ is a chain homotopy
from $\phi^j$ to $|X|\phi^{j-1}$.
\end{prop}

\begin{proof} We need to show that \[\d_{n+1}D_n^j(\tuple{x}) + 
 D_{n-1}^j\d_n(\tuple{x}) = \pm (\phi_n^j(\tuple{x}) - 
|X|\phi_n^{j-1}(\tuple{x})).\]

For $j>n$, we have $D_n^j = D^j_{n+1} = 0$, while 
\[ \phi^j_n(\tuple{x}) = |X|^{j-n}\phi_n^n(\tuple{x})=
 |X|(|X|^{(j-1)-n}\phi_n^n(\tuple{x})) =|X|\phi_n^{j-1}(\tuple{x})\]
as required.

For $j=n$, we have $D_{n-1}^j=0$ and
\begin{eqnarray+} 
 \d_{n+1}D_n^j(\tuple{x}) & = & \sum_{i=1}^{n+1} (-1)^i 
 (\d_i^0D_n^n(\tuple{x}) - \d_i^1D_n^n(\tuple{x})) \\
 & = & \sum_{i\le n} (-1)^i (\d_i^0D_n^n(\tuple{x}) - 
 \d_i^1D_n^n(\tuple{x})) \\
 & & + (-1)^{n+1}(\d_{n+1}^0D_n^n(\tuple{x}) 
 -\d_{n+1}^1D_n^n(\tuple{x})). 
\end{eqnarray+}
By Lemma \ref{lem:dDlow}, the first sum adds to zero, 
and by Lemmas \ref{lem:dDmid0} and \ref{lem:dDmid1} we have
\[ \d_{n+1}D_n^j(\tuple{x})= (-1)^{n+1}( |X|\phi_n^{n-1}(\tuple{x}) - 
 \phi_n^n (\tuple{x})),\]
as required.

For $j<n$,
\begin{eqnarray+} \d_{n+1}D_n^j(\tuple{x}) & = & \sum_{i=1}^{n+1} (-1)^i 
 (\d_i^0D_n^j(\tuple{x}) - \d_i^1D_n^j(\tuple{x})) \\
 & = & \sum_{i\le n} (-1)^i (\d_i^0D_n^j(\tuple{x}) - 
 \d_i^1D_n^j(\tuple{x})) \\
 & & + (-1)^{j+1}(\d_{n+1}^0D_n^j(\tuple{x}) 
 -\d_{n+1}^1D_n^j(\tuple{x})) \\
 & & + \sum_{i=j+2}^{n+1}(-i)^i (\d_i^0D_n^j(\tuple{x}) 
 + \d_i^1D_n^j(\tuple{x})) 
\end{eqnarray+}
which, by Lemmas \ref{lem:dDlow}, \ref{lem:dDmid0} and
\ref{lem:dDmid1} as above yields
\[\d_{n+1}D_n^j(\tuple{x})= (-1)^{j+1}( |X|\phi_n^{j-1}(\tuple{x}) - 
 \phi_n^j (\tuple{x})) + \sum_{i=j+2}^{n+1}(-i)^i (\d_i^0D_n^j(\tuple{x}) 
 + \d_i^1D_n^j(\tuple{x})) \] 
Now,
\begin{eqnarray+} 
 D_{n-1}^j\d_n(\tuple{x}) & = & 
 \sum_{i=1}^n(-1)^i (D_{n-1}^j\d_i^0(\tuple{x})  
 - D_{n-1}^j\d_i^1(\tuple{x})) \\
 & = & \sum_{i=1}^j (-1)^i(D_{n-1}^j\d_i^0(\tuple{x}) - 
 D_{n-1}^j\d_i^1(\tuple{x})) \\
 & & + \sum_{i=j+1}^n (-1)^i(D_{n-1}^j\d_i^0
 (\tuple{x}) - D_{n-1}^j\d_i^1(\tuple{x})). 
\end{eqnarray+}
The first sum is zero by Lemma \ref{lem:Ddlow}, 
and applying Lemmas \ref{lem:dDhigh0} and \ref{lem:dDhigh1} we get
\[ D_{n-1}^j\d_n(\tuple{x}) = \sum_{i=j+1}^n (-1)^i 
 (\d_{i+1}^0D_n^j(\tuple{x}) - \d_{i+1}^1D_n^j(\tuple{x})). \]
Reindexing this sum by replacing $i+1$ with $i'$, we have 
\[ D_{n-1}^j\d_n(\tuple{x}) = \sum_{i'=j+2}^{n+1} (-1)^{i'+1} 
 (\d_{i'}^0D_n^j(\tuple{x}) - \d_{i'}^1D_n^j(\tuple{x})), \]
so that 
\[ \d_{n+1}D_n^j(\tuple{x}) + D_{n-1}^j\d_n(\tuple{x}) = 
(-1)^{j+1} (|X|\phi_n^{j-1}(\tuple{x}) - \phi_n^j(\tuple{x})) \]
as required. 
\end{proof}

\begin{proof}[Proof of Theorem \ref{thm:main}] We deal first with the
case $W=R$.  There is a chain map $\maps{\pi^R}{C_*(X)}{C_*(\orb X)}$
induced by the projection of $X$ onto its orbit rack.  In Lemma 4.2 of
\cite{CJKSb}, it is proved that for 
$\boldsymbol \omega = (\omega_1,\ldots,\omega_n) \in \orb X^n$, 
the element $\sum_{z_j\in\omega_j, j = 1,\ldots, n} (z_1,\dots,z_n)$ 
of $C^R_n(X)$ is a cycle. Since the boundary maps in 
$C_*(\orb X)$ are all zero, this means that we can define a chain map 
$\maps{\psi}{C^R_*(\orb X)}{C^R_*(X)}$ by setting
\[\psi_n(\boldsymbol\omega) = \biggl(\prod_{i=1}^{n}N_{\omega_i}\biggr)
\sum_{z_j\in\omega_j, j = 1,\ldots, n} (z_1,\dots,z_n).\]
This is almost the same as the chain map used in the proof of Theorem
4.1 of \cite{CJKSb}.) Then, for $\tuple x \in X^n$,
\begin{eqnarray+}
\psi_n\pi^R_n(\tuple x) & = &\biggl(\prod_{i=1}^{n}N_{\pi(x_i)}\biggr)
\sum_{z_j\in\pi(x_j), j = 1,\ldots, n} (z_1,\dots,z_n)\\
& = & \sum_{\tuple y \in X^n} (x_1^{y_1},\dots,x_n^{y_n})\\
& = & \phi_n^n(\tuple x).
\end{eqnarray+}

Hence, by Proposition \ref{prop:htic}, the induced map
$\maps{\psi_*\pi^R_*}{H^R_n(X)}{H^R_n(X)}$ is multiplication by
$|X|^n$.
It follows, since $H^R_n(\orb X)$ is free abelian, that the torsion
subgroup of $H^R_n(X)$ is equal to $\ker \pi^R_*$ and is annihilated
by $|X|^n$, and that $\beta_n^R(X) \le \beta_n^R(\orb X)$. Since the
reverse inequality was proved in \cite{CJKSb}, the proof in the case
of rack homology is complete.

When $X$ is a quandle, the other two cases follow from the case just
proved, Theorem \ref{thm:split}, and Theorem 4.1 of \cite{CJKSb}.
\end{proof}

\section{Computations}\label{sec:calc}

In \cite{CJKSc} (Table 1), the cohomology groups 
$H^n_Q(X;\z_p)$ of some Alexander racks are given for $n=2$ or $3$
and the first few primes $p$. These racks are of the form
$\Lambda_n/(h)$, 
and the number $m$ of orbits is easily computed from Proposition 
\ref{prop:alex}(b).
For $X = \Lambda_3/(t^2+t+1)$, $m = 3$, so according to Theorem \ref{thm:main},
the dimension of
$H^2_Q(X;\z_p)$ should be 6 for $p \neq 3$, while the value in
\cite{CJKSc} is 3 in these cases. 
This led the first author to write
a C program to check the computations. Apart from $\Lambda_3/(t^2+t+1)$,
where the recomputation gave the same values as for $\Lambda_9/(t-4)$,
the results agreed with one exception, for 
$X = \Lambda_3/(t^2-t+1)$. Here \cite{CJKSc} has 
$\dim H^2_Q(X;\z_3)=0$, while the recomputation yields
$\dim H^2_Q(X;\z_3)=1$. The value 1 is in agreement with Corollary 2.4
of \cite{M}. It turns out that the disagreement is due to typographical
errors in \cite{CJKSc}, and the values just given are the ones computed by 
Carter et al.

A variant of this program computes the integral homology
of racks; we present in Table 1 the results of some calculations.
In view of Theorem \ref{thm:h23}, we give only the quandle homology, though the
program has been run to compute rack homology with the results
expected from Theorem \ref{thm:h23}.
As in \cite{CJKSc}, the racks considered are non-trivial, of order at
most 9, and of the form $\Lambda_n/(h)$ where $h$ is a monic polynomial
whose constant term is a unit in $\z_n$. The list of racks is different
from that in \cite{CJKSc} in two ways. First, we have included
$\Lambda_3/(t^2-t-1)$ and $\Lambda_2/(t^3+t^2+t+1)$. Second, it turns
out that $R_4 \simeq \Lambda_2/(t^2+1)$, 
$\Lambda_9/(t-4) \simeq \Lambda_9/(t-7) \simeq \Lambda_3/(t^2+t+1)$,
and $R_8\simeq \Lambda_8/(t-3)$ (where $\simeq$ denotes rack-isomorphism),
and we have omitted all but the first of each isomorphism class. That
$R_4 \simeq \Lambda_2/(t^2+1)$ and
$\Lambda_9/(t-4) \simeq \Lambda_9/(t-7)$ was noted in \cite{CJKSc}.
The other isomorphisms were discovered by a brute-force computation, and
that remains the only assurance we have that the racks we have listed are 
all distinct. The existence of all these isomorphisms follows from 
the next two propositions.

\begin{table}[tb]
\makebox[\textwidth]{
\begin{tabular}{|l|cc|} 
\hline
\multicolumn{1}{|c|}
{$X$\rule{0pt}{11pt}} & $H^Q_2(X)$ & $ H^Q_3(X)$ \\[2pt]
\hline
$R_3$    & 0\rule{0pt}{11pt}  & $\z_3$                            \\
$R_4$    & $\z^2\oplus\z_2^2$ & $\z^2\oplus\z_2^4$                \\
$R_5$    & 0                  & $\z_5$                            \\
$R_6$    & $\z^2$             & $\z^2\oplus\z_3^2$                \\
$R_7$    & 0                  & $\z_7$                            \\
$R_8$    & $\z^2\oplus\z_2^2$ & $\z^2\oplus \z_2^2\oplus \z_8^2$  \\ 
$R_9$    & 0                  & $\z_9$                         \\[4pt]
$\Lambda_5/(t-2)$ & 0                  & 0                        \\
$\Lambda_5/(t-3)$ & 0                  & 0                        \\
$\Lambda_7/(t-2)$ & 0                  & 0                        \\
$\Lambda_7/(t-3)$ & 0                  & 0                        \\
$\Lambda_7/(t-4)$ & 0                  & 0                        \\
$\Lambda_7/(t-5)$ & 0                  & 0                        \\
$\Lambda_8/(t-5)$ & $\z^{12}\oplus\z_2^4$ 
                                       & $\z^{36}\oplus\z_2^{24}$ \\
$\Lambda_9/(t-2)$ & 0                  & $\z_3$                   \\
$\Lambda_9/(t-4)$ & $\z^6\oplus\z_3^3$ & $\z^{12}\oplus\z_3^{12}$ \\
$\Lambda_9/(t-5)$ & 0                  & $\z_3$               \\[4pt]
$\Lambda_2/(t^2+t+1)$ & $\z_2$ & $\z_2\oplus\z_4$      \\ 
$\Lambda_3/(t^2+1)$   & $\z_3$ & $\z_3^3$              \\
$\Lambda_3/(t^2-1)$   & $\z^6$ & $\z^{12}\oplus\z_3^3$ \\
$\Lambda_3/(t^2-t+1)$ & $\z_3$ & $\z_3\oplus\z_9$      \\
$\Lambda_3/(t^2+t-1)$ & 0      & 0                     \\
$\Lambda_3/(t^2-t-1)$ & 0      & 0                     \\[4pt]
$\Lambda_2/(t^3+1)$       & $\z^2\oplus\z_2^2$ 
                              & $\z^2\oplus\z_2^6\oplus\z_4^2$ \\
$\Lambda_2/(t^3+t^2+1)$   & 0                  & $\z_2$        \\
$\Lambda_2/(t^3+t+1)$     & 0                  & $\z_2$        \\
$\Lambda_2/(t^3+t^2+t+1)$ & $\z^2\oplus\z_2^4$ 
                              & $\z^2\oplus\z_2^8\oplus\z_8^2$ \\[2pt]
\hline
\end{tabular}
} 
\caption{Some quandle homology groups.}
\end{table}

\begin{prop} If $k$ is coprime to $n$ then
$\Lambda_{n^2}/(t-(kn+1)) \simeq \Lambda_n/((t-1)^2)$.
\end{prop}

\begin{proof} We identify $\Lambda_{n^2}/(t-(kn+1))$ with $\z_{n^2}$
under the operation $a^b = (kn+1)a - knb$. There is a short exact sequence
of abelian groups
\[
 0 \to \z_n \mapright\alpha \z_{n^2} \mapright\beta \z_n \to 0,
\]
where $\alpha(1) = n$ and $\beta(1)=1$. Note that for $a\in \z_{n^2}$,
$\alpha^{-1}(na) = \beta(a)$. Let $\maps\gamma{\z_n}{\z_{n^2}}$ be a function
such that $\beta\gamma = \id$, and define $\maps\delta{\z_{n^2}}{\z_n}$
by $\delta(a) = \alpha^{-1}(a - \gamma\beta(a))$. 
The function $\z_{n^2} \to \z_n^2$ sending $a$ to $(\beta(a),\delta(a))$
is a bijection. Now define
$\maps f {\z_{n^2}}{\Lambda_n/((t-1)^2)}$ by 
$f(a) = k\beta(a) + (t-1)\delta(a)$; because $k$ is coprime to $n$,
$f$ is also a bijection. We have, for $a, b\in \z_{n^2}$, 
$\beta(a^b) = \beta(a)$ and
\begin{eqnarray*}
 \delta(a^b) & =& \alpha^{-1}((kn+1)a-knb - \gamma\beta(a)) \\
 & = & \alpha^{-1}(kna)-\alpha^{-1}(knb) + \alpha^{-1}(a - \gamma\beta(a)) \\
 & = & k\beta(a)-k\beta(b) + \delta(a).
\end{eqnarray*}
Hence
\begin{eqnarray*}
 f(a^b) & = & k\beta(a) +(t-1)(k\beta(a)-k\beta(b) + \delta(a)) \\
 & = & kt\beta(a)+(t-1)\delta(a)+k(1-t)\beta(b)
\end{eqnarray*}
On the other hand,
\begin{eqnarray*}
 f(a)^{f(b)} & = & t(k\beta(a) + (t-1)\delta(a)) + 
  (1-t)(k\beta(b)+(t-1)\delta(b)) \\
 & = & kt\beta(a) +(t-1)\delta(a) +k(1-t)\beta(b),
\end{eqnarray*}
so $f$ is the desired isomorphism.
\end{proof}

\begin{prop} If $n$ is divisible by $4$ then
$R_{2n} \simeq \Lambda_{2n}/(t - (n-1))$.
\end{prop}
\begin{proof} Here the underlying sets of both racks are naturally identified
with $\z_{2n}$. We use $a^b$ for the rack operation in $R_{2n}$,
and $a^{[b]}$ for that in $\Lambda_{2n}/(t - (n-1))$.
Thus, for $a, b \in \z_{2n}$,
\begin{eqnarray*}
 a^b & = & 2b - a \\
 \hbox{and}\qquad a^{[b]} & =& (n-1)a+(2-n)b.
\end{eqnarray*}
Define functions $\epsilon$ and $f$ from $\z_{2n}$ to itself by
\[
 \epsilon(a) = \cases{
 0,& if $a \equiv 0$ or $1 \pmod 4$;\cr
 n,& if $a \equiv 2$ or $3 \pmod 4$;\cr}
\]
and $f(a) = a +\epsilon(a)$. Since $f(a) \equiv a \pmod 4$, $f$ is
an involution. Since $a^b \equiv a \pmod 2$, we have that
$\epsilon(a^b) = \epsilon(a)$ iff $a^b \equiv a \pmod 4$,
which in turn is equivalent to $a \equiv b \pmod 2$. Since $\epsilon$
only takes on the values 0 and $n$, this implies that
$\epsilon(a^b) = \epsilon(a) + n(a-b)$. Hence
\begin{eqnarray*}
 f(a^b) & = & 2b - a +\epsilon(a) + n(a-b) \\ 
 & = & (n-1)a + \epsilon(a) +(2-n)b.
\end{eqnarray*}
On the other hand,
\begin{eqnarray*}
 f(a)^{[f(b)]} & = & (n-1)(a + \epsilon(a)) + (2-n)(b+\epsilon(b)) \\
 & = & (n-1)a + \epsilon(a) + (2-n)b,
\end{eqnarray*}
so $f$ is the desired isomorphism. 
\end{proof}

\end{document}